\documentclass[10pt]{article}

\usepackage{graphicx}
\usepackage{graphicx,psfrag}

\def \vs {\vskip 0.3cm}

\def \n {\noindent}

\def \E {{\bf  \, E\,}}

\def \T {{\hbox{Tr\,}}}

\def \CC {{\cal C}}

\def \CT {{\cal T}}

\def \CW {{\cal W}}

\def \CI {{\cal I}}
\def\CL {{\cal L}}
\def \CQ {{\cal Q}}

\def \CF {{\cal F}}

\def \tU {\tilde U}

\def \a {\alpha}

\def \b {\beta}
\def \g {\gamma}

\def \s {\sigma}

\def \d {\delta}

\def \vep {\varepsilon}
\def \e {\epsilon}

\def \z {\zeta}

\begin{document}

\title{ A Class of Even Walks and  Divergence of High Moments
of  Large Wigner Random Matrices\footnote{ {\bf Acknowledgements:} The financial support of the research grant ANR-08-BLAN-0311-11 "Grandes Matrices
Al\'eatoires" (France)  is gratefully acknowledged}\ \ 
}

\author{O. Khorunzhiy\\ Universit\'e de Versailles - Saint-Quentin, Versailles\\ FRANCE\\
{\it e-mail:} khorunjy@math.uvsq.fr}

\maketitle

\begin{abstract}
We consider the Wigner ensemble of  $n\times n$  random matrices $\hat A^{(n)}$ with truncated
to the \mbox{interval} $(-U_n,U_n)$ elements
and study
the   moments { $M_{2s}^{(n)} = \E \T(\hat A^{(n)})^{2s}$ }
by using their representation as the sums over the set of  weighted even closed walks $\CW_{2s}$.

We construct a subset $\CW'_{2s}\subset \CW_{2s}$ such that
the corresponding sub-sum diverges in the limit
$n,s\to\infty$, $s = s_n= O (n^{2/3})$ for any truncation of the form
$U_n = 2n^{1/6+\e}$ with $\e>0$, provided the probability distribution
of the matrix elements $a_{ij}$ belongs to a class of  distributions
 with the twelfth moment unbounded.
This result allows us to put forward a conjecture that the existence of  $\E \vert a_{ij}\vert^{12}$ 
is \mbox{necessary} for the existence of the
universal
upper bound of the sequence $M_{2s_n}^{(n)}$ as $n\to\infty$ and eventually for the edge spectral universality
in Wigner ensembles or random matrices. 
\end{abstract}
\vs

{\it Running title:} Even Walks and Divergence of Moments 

\vs

\section{Wigner ensembles of random matrices}
The Wigner ensemble of random matrices is given by
the family of  real symmetric (or Hermitian) random matrices
$\{A^{(n)}\}$ with  the matrix elements
$$
A^{(n)}_{ij} = {1\over \sqrt n} a_{ij},\quad i,j=1,\dots,n,
\eqno (1.1)
$$
where $\{a_{ij}, 1\le i\le j\}$ are real jointly independent
random variables of the same law such that
$$
\E a_{ij}=0, \quad \hbox{and}\quad  \E \vert a_{ij}\vert^2 = v^2.
\eqno (1.2)
$$
We denote $V_{2k}= \E  \vert a_{ij}\vert^{2k}$ for $k\ge 2$ and assume that
the probability distribution of random variables $a_{ij}$ is  symmetric.

We study the moments of the Wigner random
matrices that are obtained from (1.1) by the standard truncation procedure.
Namely, let us  consider the truncated variables
$$
\hat a_{ij}^{(n)} = \cases{ a_{ij}, & if $\vert a_{ij}\vert
< U_n$, \cr
0, & if $\vert a_{ij}\vert \ge U_n$ \cr}
$$
and determine matrices
$$
\hat A^{(n)}_{ij} = {1\over \sqrt n} \hat a_{ij}^{(n)}, \quad i,j=1,\dots,n.
$$
We will refer to the family $\{\hat A_n\}$ as to the ensemble of truncated Wigner  random matrices.
In the case of Hermitian matrices we assume that $a_{ij} = b^{(1)}_{ij}+ {\hbox{i}} b_{ij}^{(2)}$,
where $\{b^{(l)}_{ij}, l=1,2\}$ are i.i.d. random variables such that (1.2) holds and consider truncation 
 of $b^{(l)}_{ij}$.

Since the pioneering works of E. Wigner \cite{W}, one extensively uses the natural representation of
  the averaged moments of (1.1)
as the sum over a family of  closed paths of $2s$ steps.
More precisely, one can write equality
$$
\hat M_{2s}^{(n)}= {1\over n^s}
\sum_{I_{2s} \in \CI_{2s}^{(n)}} \, \Pi (I_{2s}),
\eqno (1.3)
$$
where $I_{2s} = (i_0, i_1,\dots, i_{2s-1}, i_0)$,
$\CI_{2s}^{(n)}$ is the set of all even closed trajectories $I_{2s}$
over the set $\{1,\dots,n\}$
and
$\Pi(I_{2s}) $ is the  weight given by the
mathematical expectation of the product of $a_{ij}$'s
that corresponds to $I_{2s}$,
$$
\Pi(I_{2s}) = \E \left( \hat a_{i_0,i_1}^{(n)} \cdots \hat a_{i_{2s-1},i_0}^{(n)}\right).
$$
The trajectory $I_{2s}$ can be viewed as  an $n$-realization of an even closed walk
\mbox{$w_{2s} \in \CW_{2s}$} such that in the corresponding
multigraph  $g_{2s} = g(w_{2s})$ each couple
of vertices is joined by even number $2l$ of edges,
$l\ge 0$ and the total number of edges is $2s$.
For the rigorous definitions of trajectories,
walks and their graphs see \cite{K}.

\vskip 0.3cm
In paper \cite{K}, it is proved that if there exists $\d_0>0$ such that 
 $\E \vert a_{ij}\vert^{12+\d_0}<\infty  $,
then there exists a truncation $U_n= n^{1/6-\vep}$, $\vep>0$  such that
$$
P\left\{ A^{(n)} \neq \hat A^{(n)} \ \ \hbox{infinitely often} \ \right\} = 0
\eqno(1.4)
$$
and that the moments (1.3) are bounded from above 
$$
\limsup_{n,s_n\to\infty} {1\over (4v^2)^{s_n}}
\hat M_{2s_n}^{(n)} \le \CL(\vartheta), \quad  s_n = \lfloor \vartheta n^{2/3}\rfloor, \  \vartheta>0  
 \eqno (1.5)
$$
where  $\CL(\vartheta) $ is a universal constant that does not depend
on the particular values of $V_{2k}, 1\le k\le 6$. 
Relation (1.4) is obtained by the standard application to random matrices  the Borel-Cantelli lemma,
and the bound (1.5) is proved with the help of a completed and improved version 
of the method developed by Ya. Sinai and A. Soshnikov \cite{SS1,SS2,S} to study the sum (1.3).

\vskip 0.3cm

In the present note we consider  the inverse situation. 
We prove that if the probability distribution of $a_{ij}$
is such that $\E\vert a_{ij}\vert^{12}$ does not exist, then for any
 truncation $U_n$ 
 needed in the standard proof of (1.4),
the moments of the corresponding truncated Wigner random matrices diverge in the 
asymptotic regime of (1.5).
To show this, we  construct the family $\CW'_{2s} \subset \CW_{2s}$
of even closed  walks $w'_{2s}$ such that the 
weighted
sum over the trajectories $I_{2s}$ from the corresponding equivalence classes $\CC(w'_{2s})$
$$
R^{(n)}_{2s}= {1\over n^s} \sum_{w'_{2s}\in \CW'_{2s}}
\ \sum_{I_{2s} \in \CC_n(w'_{2s})} \, \Pi(I_{2s})
\eqno (1.6)
$$
diverges. This  outcome means that the condition $\E\vert a_{ij}\vert^{12}<+\infty $
is hard to be avoided in the simultaneous proof of (1.4) and (1.5) in the frameworks of the approach  of \cite{K}.

\vskip 0.3cm

Being focused on the combinatorial part of the construction of $\CW_{2s}$,
we do not pay much attention to the generality of our results with respect to the
probability distribution $F(x)$ of $a_{ij}$ and consider one of the simplest cases
given by the probability density  
$$
F'(x) = \cases{ {C_\varphi} \,  \vert x\vert^{-13} \, \varphi(x) , & if $\vert x\vert >1$\cr
0, & otherwise,\cr} 
\eqno (1.7)
$$
where $C_\varphi$ is the normalization constant and $\varphi$ is a bounded even  positive monotone decreasing function such that 
$(\ln x)^{-1}\le \varphi(x) \le 1$ for all $x>x_0\ge 1$.

\vskip 0.1cm
Our main result  is as follows.
\vskip 0.1cm
{\bf Theorem 1.} {\it Let the probability distribution of $a_{ij}$ (1.2) has a symmetric density $f(x)=F'(x)$
of the form (1.7).
Then
for any choice of $\e>0$ in  $U_n=2n^{1/6+\e}$,
the high moments of the corresponding Wigner ensemble of truncated matrices
diverge; namely,
$$
\limsup_{n\to\infty} {1\over (4v^2)^{s_n}}
\hat M_{2s_n}^{(n)} = +\infty,
\eqno(1.8)
$$
where $ s_n = \lfloor \vartheta n^{2/3}\rfloor $ with any  $\vartheta>0$.
}

\vskip 0.1cm

Let us complete this section by the following comments.
By itself,  \mbox{Theorem 1}  does not imply  the necessity of the twelfth moment for the
bound (1.5). However, one can see that  in many cases
the truncation $\hat a_{ij}^{(n)}$ combined with the reasoning based on the Borel-Cantelli lemma  
leads to  the optimal conditions for the moments of $a_{ij}$
(see the papers \cite{G} and \cite{P} for the necessary and sufficient conditions for the semicircle law
to be valid; see \cite{BY} for the convergence of the spectral norm).

Therefore in view of Theorem 1, it  is natural to put forward a hypothesis that the 
bound $\E\vert a_{ij}\vert ^{12}<+\infty $ represents the necessary condition 
for the existence of the upper bound (1.5).
At the end of the present paper, we  discuss   the optimal conditions on $a_{ij}$
in asymptotic regimes different from that of (1.5).

Finally,  let us note that the bound (1.5) and its proof play the central role in the demonstration of the
universality of the edge spectral distribution of large Wigner random matrices \cite{S}.
Thus one can relate our results with the necessary and sufficient conditions of the
spectral universality in the spectral theory of random matrices.

\section{Proof of Theorem 1}

Let us briefly explain the  structure of the elements  $w'_{2s} \in \CW'_{2s}$ of (1.6)
which is based on a fairly  simple principle.
We assume that each walk $w'_{2s}$ consists of two parts; during  the first half-part
it creates a vertex $\z_0$ of high self-intersection degree, where $D$ edges of multiplicity $10$
meet each other. On its second part, the walk $w'_{2s}$ enters $\z_0$ and
then performs  a sequence of "there-and-back"
trips along the multiple edges. The vanishingly small proportion of the
half-part trajectories with $\z_0$ of this kind is compensated 
by the large number of choices
where to go by these
there-and-back
steps form the second half-part. The
cardinality of trajectories of this structure is such that it is  impossible
to get a finite upper bound  for the sub-sum (1.6)
with $U_n= 2n^{1/6+\e}$, $\e>0$.


\subsection{Construction of the walks $\CW'_{2s}$}

Let us split the time interval $[0,2s]$
into two  parts,
the $X$-part containing  $2s'+2$ steps and
the $Y$-part containing $2L$ steps
with obvious equality \mbox{$2s'+2 +2L = 2s$.}
\begin{figure}[htbp]
\centerline{\includegraphics[width=11.5cm,height=3.5cm]{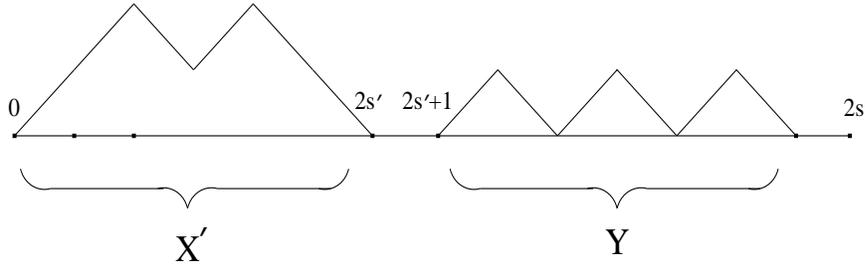}}
\caption{\footnotesize{ Partition of $[0,2s]$ into $X'$- and $Y$-parts and Catalan structures }}
\end{figure}
The chronological order is as follows:
the first $2s'+1$ steps belong to the $X$-part;
these are followed by $2L$ steps of the $Y$-part
and the final step  $[2s-1,2s]$ is attributed again to
the $X$-part.
We determine in  a natural way
the  sub-walks that correspond to the $X$- and $Y$-parts of $w_{2s}$
and denote them by $w^{(X)}_{2s}$ and $w_{2s}^{(Y)}$, respectively.
We also  consider
the first $2s'$ steps as the $X'$-part of the walk  (see figure 1).

\subsubsection{Catalan structures on $X'$-part}

Regarding the $X'$-part, one can  create a Catalan structure by   attributing  to the $2s'$ steps $s'$ signs
$"+"$ and $s'$ signs $"-"$ \cite{St}. Following \cite{SS1}, we will say that
each step
labelled  by $+$ is marked,
the remaining steps being the non-marked ones.
The Catalan structure is in one-to-one correspondence
with the set $\CT_{s'}$ of plane rooted trees  $T$
constructed with the help of $s'$ edges
or equivalently, with the set of the Dyck paths \cite{St}.
This correspondence can be determined
by the lexicographic (or in other words,
by the chronological) run over the element  $T\in \CT_{s'}$
that produces also a  Dyck path.
In this representation, the signs $+$ and $-$ correspond to the ascending and the descending steps of the Dyck path,
respectively.

It is known that the cardinality of $\CT_{s'}$
is given by the Catalan number,
$$
\vert \CT_{s'}\vert  = C^{(s')} = {(2s')\over s'! (s'+1)!}.
\eqno (2.1)
$$
In our construction, we are mainly related with  a subset $\tilde \CT_{s'}$ of Catalan structures 
on $X'$-part
 such that the majority of  vertices of the corresponding trees have
the degree $6$ with $5$ children edges.   To describe  this subset,
 we introduce the number
$s'' = \lfloor s'/5\rfloor$
and consider the set of all Catalan trees $\CT_{s''}$.
Given $T\in \CT_{s''}$, we start the lexicographic run
over $T$
and after each descending step $-$ we add
the sequence $(+,-,+,-,+,-,+,-)$.
Then we get  a Dyck path of $10s''$ steps.
We denote by $T_4\{\CT_{s''}\}$
the set of all Dyck paths obtained by this procedure.
To get the subset $\tilde \CT_{s'}$, we
consider the elements of  $T_4\{\CT_{s''}\}$
added by all possible Catalan structures on the remaining $2s'-10s''$ steps.

The final specification is that we are going to use  the subset  of trees $\CT^{(d_0)}_{s''} \subset \CT_{s''}$ that
have vertices with the number of children not greater than a given value  $d_0>2$.
The following  estimate from below
$$
\vert \CT^{(d_0)}_{s''} \vert \ge \left(1 - (2s''+1) e^{-(d_0-2) \ln (4/3)}\right) C^{(s'')}
\eqno (2.2)
$$
can be proved, in particular, by using the recurrence relations for $C^{(s)}$ (see \cite{K}).

The subset we need $\tilde \CT^{(d_0)}_{s'}$ is constructed with 
the elements of $T_4\{\CT^{(d_0)}_{s''}\}$ completed on
the remaining $2s'-10s''$ steps by all possible Catalan structures.
We will say that the edges of $T= \tilde \CT^{(d_0)}_{s''}$
represent the principal edges of the corresponding element of
$\tilde \CT^{(d_0)}_{s'}$, while the others are the supplementary ones.

Using (2.1) with $s''$ and taking into account  inequalities
$$
\sqrt{2\pi k}   \left({k\over e}\right)^k \le k! \le e\sqrt{2\pi k}     \left({k\over e}\right)^k,
\eqno (2.3)
$$
we can write that in the limit of large $s'$,
$$
\vert \tilde \CT_{s'}\vert \ge {(2s'')!\over s''! (s''+1)!}
\ge 
{4^{s'/5 +3}\over 2e^2\sqrt \pi }
\cdot {1\over (s'+1)^{3/2}}.
$$
Then we can deduce from (2.2) the estimate 
$$
\vert \tilde \CT^{(d_0)}_{s'} \vert 
\ge  \left(1 - s'' e^{-(d_0-2) \ln (4/3)}\right) 
{4^{s'/5 +3}\over 2e^2\sqrt \pi }
\cdot {1\over (s'+1)^{3/2}}.
\eqno(2.4)
$$

\subsubsection{Self-intersections on $X'$-part}

Regarding the $X$-part of $w_{2s}$, 
we will say that the graph $\tilde g(w^{(X)}_{2s})$
with simple non-oriented edges
represents the {\it frame} of the walk $w_{2s}$
that we  denote by $\CF(w_{2s})$.
Using the Catalan structure on the $X'$-part,
we construct  $w^{(X)}_{2s}$ in the following way.

At the zero instant of time, we put the walk at the root vertex $w(0)= \rho$.
Here and below we simply say that $\{w_{2s}(t), 0\le t\le 2s\}$ represent the vertices of the graph $g(w_{2s})$
instead of more rigorous but   cumbersome
expression.

We start the run over  the Catalan structure $T$; if the step $[t,t+1]$
is marked, then  the walk either  creates a new vertex $w(t+1)$
or arrives at  one of the already existent vertices. We say that  the edge
$(w(t), w(t+1))$ is marked.
If the step $[t,t+1]$ is such that $w(t+1)$  joins an existing vertex $\b$,
then we say that $t+1$ is the instant of the {\it self-intersection} of the walk \cite{SS1}.
The total number of arrivals at $\b$ at the marked instants of time
is called the {\it self-intersection degree} of $\b$; we denote this number by $\kappa(\b)$.
If $\kappa(\b)=2$, then we say that $\b$ is the vertex of {\it simple self-intersection}.

If the instant of time $[t,t+1]$ is non-marked,
then the walk performs the step along the marked edge $(\b,\g) $ attached to $w(t)=\b$ 
such that  $(\b,\g)$ is passed odd number of times during the interval $[0,t]$.
This marked edge is uniquely determined. 

\vskip 0.1cm
Let us consider a particular tree $T_{s'}\in \tilde \CT^{(d_0)}_{s'}$
together with  its predecessor \mbox{$T_{s''}\in \CT_{s''}$} and choose
$D$ edges $e'_1, \dots e'_D$
among the principal edges of this $T_{s'}$.
We want these edges to be such that for any couple
$e'_i$ and $e'_j$,
the  distance between their vertices in $T_{s''}$ is not
less than 3.
It is easy to see that the number of possible choices is estimated from below by
$$
{1\over D!} s''\left(s''-2d_0^3\right)
\left(s''-3d_0^3\right) \cdots
\left(s''-(D-1)d_0^3\right)\ge {1\over D!} \left( s''- Dd_0^3\right)^D.
\eqno (2.5)
$$

Let us denote by $\tau_1< \dots < \tau_D$ the instants of time that correspond
to the edges $e'_j$, $1\le j\le n$  in $T_{s'}$.
Let us denote by $\tilde \tau^{(j)}_i$, $ 1\le i\le 4$ the marked instants of time that correspond
to the supplementary edges of $T_{s'}$ that follow after that  $e'_j$
is passed for the second time.  With particular values of $\tau_j$ pointed out,  
we force the walk
to join  the same vertex $\z_0 = w(\tau_1)$ at the instants of time
$\tau_j, 2\le j\le D$.
Also we oblige the walk
to join  $\z_0$ at the instants of time 
$\tilde \tau^{(j)}_i$.
Then we get a walk that has a vertex 
of the self-intersection degree
$\kappa(\z_0)=5D$; there are $D$ distinct vertices $\xi_j$
such that the edge $(\xi_j,\z_0)$ is passed
10 times by $w_{2s}$ when counted in both directions.

\vskip 0.3cm
The last stage is to create $\nu_2$ simple self-intersections
with the help of 
 $s' - 5D$ marked steps not used before.
This can be done in not more than 
$$
{1\over \nu_2!} {(s'-5D)(s'-5D-1)\over 2 }
\cdots {(s'-5D-2\nu_2+1)(s'-5D-2\nu_2)\over 2}
$$
$$
\ge {1\over \nu_2!}
\left( {(s'-5D-2\nu_2)^2\over 2}\right)^{\nu_2}
\eqno (2.6)
$$
ways.

The $X'$-part of the walk being constructed,
we attribute the sing $+$ to the step $[2s',2s'+1]$;
the prescription is such that $w_{2s}(2s'+1) = \z_0$.
It is clear that $w(2s') = \rho$
and therefore the  step $[2s', 2s'+1]$
produces the edge $(\rho,\z_0)$.
The last step of the $X$-part $[2s-1,2s]$  is non-marked; it returns the walk from $\z_0$ to
the root vertex $\rho$.

\subsubsection{The core of the walk and $w_{2s}^{(Y)}$}

\begin{figure}[htbp]
\centerline{\includegraphics[width=12cm,height=4cm]{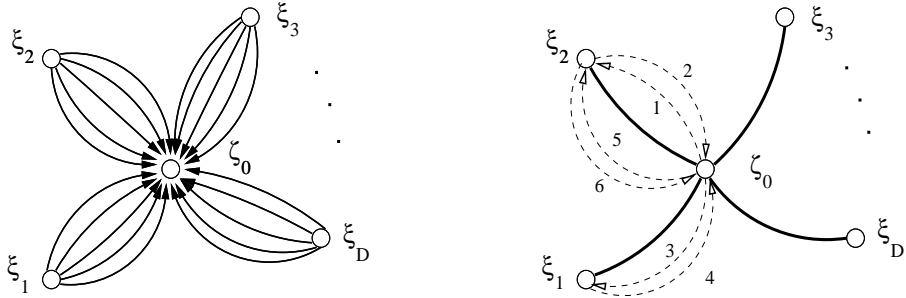}}
\caption{\footnotesize{ The multiple edges, the core $\CQ^{(5)}$ and 6 first steps of the $Y'$-part over the core}}
\end{figure}
We see that the walk $w_{2s}$ at the instant
$t=2s'+1$ enters the vertex $\z_0$ that belongs to the
family of multiple edges   $\CQ=\{(\xi_i,\z_0), \ i=1,\dots, D\}$, where  
$\xi_i\neq \xi_j$,
each edge being passed by the walk $2p=10$ times.

We refer to the   set of corresponding simple non-oriented edges  
as to the {\it core }
of the walk and denote it by $\CQ^{(p)} = \CQ^{(5)}$.
The $Y$-part $w^{(Y)}_{2s}$ 
represents the sequence of $L$ "there-and-back"
trips of the form $(\z_0,\xi_j,\z_0)$ over the core (see figure 2).
\vskip 0.3cm

To get the appropriate estimates for the weights,
we obligate the first  $2D$ steps of $w_{2s}^{(Y)}$ to  perform the ordered  trips $(\z_0,\xi_1,\z_0)$,
$(\z_0,\xi_2,\z_0)$, $\dots (\z_0,\xi_D,\z_0)$.
As for the remaining $2L-2D$ steps that determine the $Y'$-part of the walk,
there is no restrictions on the choice
to which of the vertices $\xi_j$ to go; the total
number of
all possible  walks $w^{(Y')}_{2L-2D}$ is obviously given by $D^{L-D}$.

\vskip 0.3cm
Regarding the combination of all possible sub-walks $w_{2s'}^{(X)}$ and $w_{2L}^{(Y)}$
and using inequalities  (2.5) and (2.6), we obtain the following estimate from below for the cardinality of $\CW'_{2s}(d_0,D,L,\nu_2)$:
$$
\vert \CW'_{2s}(d_0,D,L,\nu_2)\vert \ge \vert \tilde \CT^{(d_0)}_{s'}\vert 
\cdot {\left( s''- Dd_0^3\right)^D\over D!} 
\cdot {1\over \nu_2!}
\left( {(s'-5D-2\nu_2)^2\over 2}\right)^{\nu_2} \cdot D^{L-D},
\eqno (2.7)
$$
where 
$ s'' = \lfloor {s'/ 5}\rfloor$.

\subsection{Trajectories and weights}

It is easy to see that the graph $g_{2s}= g(w_{2s})$
of the walk $w_{2s} \in \CW'_{2s}(d_0,D,L,\nu_2)$
we constructed has exactly $s+1 - 5D - \nu_2-L $ ordered vertices.
Therefore given   $w_{2s} \in \CW'_{2s}$ and assigning
different values from the set $\{1,\dots, n\}$
to the vertices of $w_{2s}$, we  get the equivalence
class of trajectories $\CC(w_{2s})$
of the cardinality
$$
\vert \CC(w_{2s})\vert =
{n(n-1) \cdots  (n-(s-5D-\nu_2-L))}.
$$
Taking into account the factor $n^{-s}$,
we can write that
$$
{1\over n^s}
\vert C(w_{2s})\vert =
{1\over n^{5D + \nu_2+L} }
\cdot \prod_{i=1} ^{s-5D-\nu_2-L}
\left( 1-{i\over n}\right)
$$
$$
\ge {1\over n^{5D + \nu_2+L} }\cdot  \exp\left\{- {(s'+2-\nu_2-5D)^2\over 2n}\right\}.
\eqno (2.8)
$$
Indeed, since we are related with the asymptotic regime
$s= O(n^{2/3}), n\to\infty$, we can write that
$$
 \prod_{i=1} ^{s-5D-\nu_2-L}
\left( 1-{i\over n}\right)=
\exp\left\{ \sum_{i=1}^{s-5D-\nu_2-L} \ln \left( 1 - {i\over n}\right)\right\}
$$
$$
\ge \exp\left\{ -\sum_{i=1}^{s-5D-\nu_2-L}{i\over n}\right\}
 \ge \exp\left\{ - { (s-5D-\nu_2 -L+1)^2\over 2n}\right\}.
$$
Then, taking into account equality $s-L = s'+1$, we get (2.8). 

\vskip 0.3cm
Regarding the weights of the trajectories, we observe
that if $I'_{2s}$ and $I''_{2s}$ belong to the same
equivalence class $\CC(w_{2s})$, then
$$
\Pi(I'_{2s}) = \Pi(I''_{2s}) = \Pi(w_{2s}).
$$
Each edge $(\a,\b)$ of the frame $\CF(w_{2s})$ represents
a random variable $\hat a_{(\a,\b)}$ and these random variables
are jointly independent. Each edge of the frame
being passed in both directions by equal number $k$
of times, we  get the factor $\E \vert \hat a_{(\a,\b)}\vert^{2k}$
with no difference between the cases of Hermitian
random matrices $\hat A^{(n)}$ and the real symmetric $\hat A^{(n)}$.

If the edge $(\a,\b)$ does not belong to $\CQ^{(5)}(w_{2s})$,
then it produces either the factor $\E \vert \hat a_{(\a,\b)}\vert^2$
or the factor $\E \vert \hat a_{(\a,\b)}\vert^4\ge 
\left( \E \vert \hat a_{(\a,\b)}\vert^2\right)^2$.
Taking into account that $\E \vert \hat a_{(\a,\b)}\vert^2 \ge v^2/2$ for large values of $n$,
we can write that
$$
\prod_{(\a,\b)\in \CF\setminus \CQ^{(5)}} \E \vert \hat a_{(\a,\b)}\vert 
^{2k_{(\a,\b)}} \ge (v^2/2)^{s-5D}.
$$

\vskip 0.3cm
Regarding the edges of the core of $w_{2s}$, we get the factor
$$
\prod_{i=1}^D \E \vert \hat a_{(\xi_i,\z_0)}\vert^{12+2l_i}, 
\quad \sum_{i=1} l_i = L-D,
$$
where $l_i\ge 0$ are determined by the walk $w_{2s}$.  
Elementary computations based on (1.7) show that
$$
\E \vert \hat a_{(\xi_i,\z_0)}\vert^{12+2l_i} \ge 2C_\varphi U_n^{2l_i} \int_{1/2}^1 {y^{2l_i-1}\over \ln y + \ln U_n} dy
\ge {2C_\varphi  \over  \ln U_n} \cdot \tU_n^{2l_i},
\eqno (2.9)
$$
where $\tU_n = U_n/2 = n^{1/6+\e}$. 

Summing up, we see that the weight of the walk $w_{2s}\in \CW'_{2s}$
is bounded from below by
$$
\Pi(w_{2s})\ge \left( {v^{2}/ 2}\right)^{s} \cdot 
\left( {2^6 C_\varphi\over v^{10} \ln U_n}\right)^D \tU^{2L-2D}.
$$
Let us note that this bound based on (2.9) could be 
relaxed; then one gets \mbox{Theorem 1} in more general situation than
that determined by (1.7). We do not discuss these generalizations here.

 Now we are ready to get the estimate from below for the sum 
$R_{2s_n}^{(n)}$ (1.6) and therefore for the moments $\hat M_{2s_n}^{(n)}$ (1.3).

\subsection{Estimate from below for $R_{2s}^{(n)}(D,L)$}

Given the particular values of $D$ and $L$ and performing the sum
over $\nu_2$ from $0$ to $\s\le s'-5D$, we can write that
$$
R_{2s}^{(n)}(D,L) = \sum_{w_{2s}\in \CW'_{2s}(d_0, D,L)} 
n^{-s} \vert \CC(w_{2s})\vert \cdot \Pi(w_{2s})
$$
$$
\ge (v^{2}/2)^{s}\  \sum_{\nu_2=0}^\s
\vert \CW'_{2s}(d_0, D,L,\nu_2) \vert 
\cdot \exp\left\{ - {(s'+2-\nu_2-5D)^2\over 2n}\right\} 
$$
$$
\times \left( {2^6C_\varphi\over v^{10} \ln U_n}\right)\cdot 
{ \tU_n^{2L-2D}\over n^{5D+L+\nu_2}}.
\eqno (2.10)
$$

Regarding the sum of (2.6) over $\nu_2$, we assume $\s = yn^{1/3}$  and  get 
with the help of (2.2) the following inequalities
$$
\sum_{\nu_2=0}^{\sigma}
{1\over \nu_2!} \left( { (s'-2\nu_2-5D)^2\over 2n}
\right)^{\nu_2}  \ge
\exp\left\{ { (s' -2\s - 5D)^2\over 2n}\right\} - \sum_{l=0}^{\infty}
{1\over l!} \left( { (s'-2\s-5D)^2\over 2n}\right)^{\s+l}
$$
$$
\ge
\exp\left\{ { (s' -2\s - 5D)^2\over 2n}\right\} \left( 1 - {1\over \sqrt{2\pi \s}}
\left( {e(s'-2\s-5D)^2\over 2n \s}\right)^{\s}\right)
$$
$$
\ge
\exp\left\{ { (s' -2\s - 5D)^2\over 2n}\right\}
\left( 1 - \left( {e \vartheta^2\over 8 y}\right)\right)\ge {1\over 2} \exp\left\{ { (s' -2\s - 5D)^2\over 2n}\right\}
$$
provided $s' = \lfloor s/2-1\rfloor  \le \vartheta n^{2/3}/2$ and $y> e \vartheta^2/4$.
Using   this inequality, we get from relations  (2.7) and (2.10) with the help of (2.3), our main bound
$$
R_{2s}^{(n)}(D,L) \ge (v^{2}/2)^s\   \vert \tilde \CT^{(d_0)}_{s'}\vert \cdot 
{1\over e\sqrt {2\pi D}}
\left( {e(s'/5-1-d_0^3D)\over Dn^6}\right)^D
$$
$$
\times\ \left( {D\tU_n^2\over n}\right)^{L-D}  \left( {2^6C_\varphi\over v^{10} \ln U_n}\right)^D .
\eqno (2.11)
$$
Let  \mbox{$L  = \lfloor \vartheta n^{2/3}/2\rfloor$}.
Given any positive $\e$ in the truncation $U_n$, we 
choose $0<\vep'<\e$ and set 
$$D= \lfloor n^{2/3-\vep'}\rfloor+1. 
\eqno (2.12)
$$ 
Using (2.4) with  $d_0 = n^{\vep'/6}$ and taking into account that 
$L-D \ge (\vartheta/3) n^{2/3}$ for large values of $n$, we deduce from (2.11) the following inequalities  
$$
{1\over (4v^2)^s} R_{2s}^{(n)}(D,L)\ge  
{1\over 8e^3\pi } \cdot
{1\over (\vartheta n^{2/3}+1)^{3/2} \, n^{1/3-\vep'/2}}
\cdot n^{\vartheta \vep'n^{2/3}/3}
$$
$$
\times \left( { 64eC_\varphi\over v^{10}( \ln n +\ln 2) \n n^6} 
\left( {\vartheta\over 10} - {1\over n^{\vep'/2}} - {2\over n^{2/3}}\right)\right)^{n^{2/3-\vep'}} 
$$
$$
\ge C_1 \exp\left\{ n^{2/3}\left( {\vartheta \vep'\over 3}\ln n -\ln 4\right) - n^{2/3-\vep'} \left( ( 6 \ln n + \ln \ln n) + C_2\right) 
- {4\over 3} \ln n
\right\},
$$
where $C_1$ and $C_2$ are the constants that depend on $v^2$, $\vartheta$ and $C_\varphi$ only.
The last expression obviously diverges as $n\to\infty$. Theorem 1 is proved.

\section{Concluding  remarks}

Let us  look at the main bound (2.11) and consider the principal  
factor 
$$ 
\left( { D U_n^2\over 4 n }\right)^{L-D}
\eqno (3.1)
$$
that causes the divergence of $R_{2s_n}^{(n)}(D,L)$ for $D= O(n^{2/3-\vep'})$ and $L= O(n^{2/3})$, $ n\to\infty$. 
We see that if one switches from the asymptotic regime $s_n = O(n^{2/3})$ to
another one given by $s_n = O(n^{\eta})$ with  $0<\eta <2/3$, then the truncation of the form 
$$
U_n = n^{1/2-\eta/2}
\eqno (3.2)
$$
 represents the critical exponent beyond that $R_{2s}^{(n)}(D,L)$ is divergent. This truncation 
guarantees (1.4) provided $\E \vert a_{ij}\vert^{4/(1-\eta)+\delta } = V_{4/(1-\eta)+\delta}$ exists with  $\delta>0$. 
The moments (1.3) can be studied in this situation by using the 
combination of the method of \cite{K} with the 
corresponding frame representation of the graphs of walks.  The divergence of the moments in the case of 
$\E \vert a_{ij}\vert^{4/(1-\eta)} = +\infty$ can be proved by construction of the 
 class of walks of the type $\CW'_{2s}$ 
with the appropriately determined  core $\CQ^{(p)}$ that gives the factors of the form (3.1).

\vskip 0.1cm 

Regarding the asymptotic regime $s_n = \lfloor \ln n\rfloor$ used in the  studies of  
the spectral norm of $A^{(n)}$, we observe that the critical value for $U_n$ (3.2) is related with 
the exponent $1/2$ that requires the existence of the fourth moment $\E \vert a_{ij}\vert^4$. 
This condition is shown to be  sufficient for the existence of 
$\limsup_{n\to\infty} \Vert A^{(n)} \Vert $  obtained  with the help of (1.4) and (1.5) with $s_n = O(\ln n)$,
and also is proved to be  the necessary one \cite{BY}. 

\vskip 0.3cm

These observations reveal once more a
 special role played by the  walks $ \CW'_{2s}$ in the 
studies of high moments of Wigner random matrices. 
The interpolation argument with respect to the  asymptotic regimes considered above
 also gives 
 more support to our conjecture that 
the existence of the twelfth moment $\E \vert a_{ij}\vert^{12}$ is the necessary and sufficient 
condition for the existence of the universal bound $\CL(\theta)$ of (1.5) as well as for the
edge spectral universality of the Wigner ensembles of random matrices.

\end{document}